\documentclass{amsart}
\usepackage{amssymb}
\usepackage{amsthm}
\usepackage {amscd}
\newtheorem{satz}{Theorem}
\newtheorem{lemma}{Lemma}

\newtheorem{kor}{Corollary}
\newtheorem{proposition}{Proposition}

\theoremstyle{remark}
\newtheorem{rem}{Remark}
\newtheorem{ex}{Example}

\newcommand{\C}{{\mathbb C}}
\newcommand{\R}{{\mathbb R}}
\newcommand{\Z}{{\mathbb Z}}

\newcommand{\CP}{\C P}
\newcommand{\cpq}{\overline{\C P^2}}

\DeclareMathOperator{\degree}{deg}

\begin{document}

\title{Immersions of surfaces in almost--complex 4--manifolds}
\author{Christian Bohr}
\address{Mathematisches Institut der LMU, Theresienstr. 39 , 80333 M\"unchen , Germany}
\email{bohr@rz.mathematik.uni-muenchen.de}
\thanks{The author has been supported by the Graduiertenkolleg ``Mathematik im Bereich
ihrer Wechselwirkung mit der Physik'' at the University of Munich}

\begin{abstract}
In this note, we investigate the relation between double points and complex points of
immersed surfaces in almost--complex 4--manifolds and show how estimates for the minimal
genus of embedded surfaces lead to inequalities between the number of double points and
the number of complex points of an immersion. We also provide a generalization of a classical
genus estimate due to V.A. Rokhlin to the case of immersed surfaces.
\end{abstract}

\maketitle

\section{Introduction}

Suppose that $X$ is a 4--manifold with an almost complex structure $J$ and $F \hookrightarrow X$ 
is an immersed oriented surface. For a generic immersion, there are two types of distinguished
points on this surface. On the one hand, we have
the singularities of the embedding which we assume to be ordinary double points (this is
always the case for a generic immersion). On the other hand, there is a finite number
of complex points, i.e. points at which the almost complex structure $J$ preserves the
tangent space of the surface. At such a point, the orientation induced by $J$ on the tangent
space may coincide with the orientation of the surface -- in which case we will call the
point a {\bf positive complex point} -- or may not, then it is called a {\bf negative complex point}.

In this paper, we use results of H.F. Lai
to derive relations between the number of double
points and the number of complex points of an immersion, thus extending the results
of \cite{CG}, where only the case of embedded surfaces is treated. Our main result concerns
immersed surfaces in the neighborhood of almost complex submanifolds.
So let us assume that $F_0 \subset X$ is an almost complex curve, i.e. an embedded oriented surface whose
points are all positive complex points.
Then clearly the number $n^-$ of negative complex points is zero, and so is the number $d_-$ of double points
having negative sign. In particular, we have the inequality $n^- \leq d_-$. It turns out that a similar
inequality is true for immersed surfaces ``near'' $F_0$.

\begin{satz}\label{diskbundle}
Let $X$ be an almost complex 4--manifold and suppose that $F_0 \subset X$ is 
an embedded pseudoholomorphic curve with $F_0 \cdot F_0 > 0$. Now let $F$
be an immersed surface contained in a tubular neighborhood
of $F_0$. Then the following inequalities hold:
\begin{enumerate}
\item If $F \cdot F_0 > 0$, then $n^-(F) \leq d_-(F)$.
\item If $F \cdot F_0 <  0$, then $n^+(F) \leq d_-(F)$.
\end{enumerate}
\end{satz}

Here $n^+$ respectively $n^-$ denote the numbers of positive and negative complex points of $F$ -- counted
properly, see Section \ref{complexpoints} for details -- and $d_{\pm}$ denotes the number of positive
respectively negative double points. The proof of this theorem, which relies on Lai's results and 
facts from Seiberg--Witten gauge theory, will be given in Section \ref{proof}. In Section \ref{complexpoints},
we review the results of Lai we will need and relate them to the existence results for pseudoholomorphic
curves proved in \cite{B1}. In Section \ref{estimates}, we generalize one of the genus estimates
given by Rokhlin in
\cite{Ro} to the case of immersed surfaces and use this to derive an inequality between the numbers
of complex points and double points of certain immersions.

Finally, I would like to thank my advisor, Dieter Kotschick, for many helpful discussions
and suggestions.

\section{Complex points of immersions}\label{complexpoints}

In this section, we will briefly describe the paper \cite{Lai} of 
H.F. Lai which contains
a formula for the algebraic number of complex points of a surface in an almost complex
4--manifold, and show how this formula is related to the results of
\cite{B1}.
We will then apply Lai's results to derive some relations between the number of complex points
and the number of double points of immersed surfaces in almost complex 4--manifolds.

Suppose that $\eta$ is a complex vector bundle of rank $n$ over an
oriented manifold $X$ and that $\eta$ splits as a direct sum
\[
\eta=\xi \oplus \xi'
\]
of a complex bundle $\xi$ of rank $(n-1)$ and a complex line bundle $\xi'$.
Then the Chern product formula for direct sums implies an obvious relation between the 
Chern classes of the bundles
$\eta,\xi$ and $\xi'$.
The question Lai examined in his paper is the following. Suppose that we have
a {\em real} subbundle $\xi \hookrightarrow \eta$ of real dimension $k=2n-2$. Then
again $\eta$ splits as above, with the important difference that the splitting is
now a splitting as a real bundle of rank $2n$. Are there still relations between
the Euler classes of the bundles $\xi$ and $\xi'$ and the Chern classes of $\eta$?

For this purpose, he uses a certain notion of ``complex point'' which we will now explain.
Consider the bundle $G_k(\eta)$ whose fibre over a point $x \in X$ is the Grassmannian
of oriented $k$--dimensional real subspaces of the fibre $\eta_x$. The inclusion
$\xi \hookrightarrow \eta$ defines a section (``Gauss map'') 
\[
t : X \rightarrow G_k(\eta)
\]
in this bundle, given by $t(x)=\xi_x$. 
Let $G^\C_{n-1}(\eta)$ denote the bundle of complex subspaces of dimension $n-1$ 
in the fibres of $\eta$. Note that every
complex subspace carries a canonical orientation and is therefore an {\em oriented}
$k$--dimensional real subspace. Hence we have a canonical inclusion $G^\C_{n-1}(\eta) \rightarrow G_k(\eta)$
whose image will be denoted by
$K^+_\eta$ (this is $K_\eta$ in \cite{Lai}). 
Since the fibres of $G^\C_{n-1}(\eta)$ are complex manifolds and $X$ is oriented,
$K^+_\eta$ carries a natural orientation and therefore 
defines a homology class in $H_*(G_k(\eta))$. In a similar manner, we can
fix an orientation on the fibre of $G_k(\eta)$ (for example given by the Schubert
calculus) to obtain an orientation of the total space $G_k(\eta)$.
If we equip every complex subspace in $\eta$ with the opposite, non--complex 
orientation, we obtain a second embedding of $G^\C_{n-1}(\eta)$ into $G_k(\eta)$
whose image will be denoted by $K^-_\eta$. Note that if 
$\nu : G_k(\eta) \rightarrow G_k(\eta)$ denotes the involution given by reversing the orientation,
$K^-_\eta = \nu (K^+_\eta)$. We will orient $K^-_\eta$ such that $\nu$ maps
the orientation of $K^+_\eta$ onto minus the orientation of $K^-_\eta$.
Lai now proved the following result (Theorem 5.10 in \cite{Lai}).

\begin{satz}[Lai]\label{laithm}
Let $\eta \rightarrow X$ be a complex vector bundle of rank $n$ and
$\eta=\xi \oplus \xi'$ a splitting (as real bundle) into oriented real vector bundles
$\xi$ and $\xi'$ of ranks $k=2n-2$ and $2$. Suppose that the orientations of $\xi$
and $\xi'$ are compatible with the complex orientation of $\eta$ 
(i.e. an oriented basis of $\xi$ together with an oriented basis of $\xi'$ defines
an oriented basis for $\eta$), and let
$t : X \rightarrow G_k(\eta)$ denote the ``Gauss section'' defined by $\xi$.
Then
\[
e(\xi) + \sum_{r=0}^{n-1} e(\xi')^r \cup c_{n-r-1}(\eta)= 2t^* PD(K^+_\eta),
\]
where $PD$ denotes Poincar\'e duality.
\end{satz}

\noindent
Note that Lai's formula also implies a statement about $t^*PD(K^-_\eta)$, namely
\begin{equation}\label{laiorg}
e(\xi)+\sum_{r=0}^{n-1}(-1)^{r+1}{e(\xi')}^r \cup c_{n-r-1}(\eta)=2t^*PD(K^-_\eta),
\end{equation}
which can easily be derived from Theorem \ref{laithm} by reversing the
orientations of $\bar{\xi}$ and $\bar{\xi'}$.

As an application of his formula, Lai considers a complex manifold $X$ and an immersed
surface $F \subset X$. He then proves an equation involving the Euler classes of the
normal bundle of $F$, the Euler class of its tangent bundle, the Chern class of $X$ and
the algebraic number of complex points of $F$. 
In \cite{CG}, Chkhenkeli and Garrity observed that Lai's arguments still hold if we consider
almost complex manifolds instead of complex manifolds, since he only deals with
vector bundles but does not make use of the fact that the complex structure on $TX$ is
integrable. However, there seems to be some confusion about the signs and the question
how to count complex points in \cite{CG}, so we work out this point in greater detail.

Let $X$ be a 4--manifold which carries an 
almost complex structure $J : TX \rightarrow TX$. We orient $X$ using the orientation
given by $J$. Assume that $\iota : F \hookrightarrow X$ is an immersion of a 
connected and oriented surface $F$ into $X$. 
If -- as above -- $G_2(TX)$
denotes the bundle of Grassmannians of 2--dimensional oriented real subspaces, the immersion
defines a Gauss map $t_F : F \rightarrow G_2(\iota^*TM)$. 

Now let us consider 
the submanifolds
$K^\pm_{\iota^*TM}$ in $G_2(\iota^*TM)$, which will be abbreviated by $K^\pm$ in the sequel.
The points $x \in F$ with $t_F(x) \in K^+$ are exactly the points where $J$ respects the tangent
space of $F$ (i.e. $J$ commutes with $d\iota$) and the orientation induced by $J$ on $T_xF$
equals the orientation of $F$. Similarly, $t_F(x) \in K^-$ means that $J$ preserves $T_xF$, but
induces the opposite orientation on $T_xF$. 
We will call the points of the first kind {\bf positive complex points} and the points of
the second kind {\bf negative complex points}.

For a generic immersion, $t_F$ will be transversal to $K^\pm$, and 
since  $K^\pm$ has codimension $2$ in $G_2(\iota^*TM)$, we have
well defined intersection numbers $n^\pm$ between $t_F(F)$ and $K^\pm$, given
by the relation
\begin{equation}\label{defnofalgsum}
n^\pm=\langle t_F^*PD([K^\pm]) , [F] \rangle.
\end{equation}
We will refer to these
numbers as the algebraic sums of positive respectively negative complex points.
Now suppose that $\iota$ is an embedding. Then we have
\[
\iota^* TM = N \oplus TF,
\] where $N$ denotes the normal bundle, and $\langle e(N),[F] \rangle = F \cdot F$ equals
the self--intersection number of $F$. Using Lai's result, applied to $\iota^*TM$, we therefore 
obtain the following equation:
\begin{equation}\label{lai}
g(F)+n^+ = 1 + \frac{1}{2} (F \cdot F + \langle \iota^* c_1(X),[F] \rangle )=
1 + \frac{1}{2} (F \cdot F - K \cdot F).
\end{equation}
Using the equation for $t_F^*K^-$ derived above, we also obtain
\begin{equation}\label{laimod}
g(F)+n^- =1 + \frac{1}{2} (F \cdot F + K \cdot F).
\end{equation}
A nice example is a holomorphic curve $F$ in a complex surface $X$. 
Then the image of the Gauss map does not meet $K^-$ at all, but is entirely contained in
$K^+$. Therefore $n^-=0$, and equation~\eqref{laimod} is just the adjunction equality. So it
seems that the ``$F \cdot C$'' used by Chkhenkeli and Garrity should be our $n^-$ to obtain
the correct results.

Using equation~\eqref{laimod}, we are now able to give another formulation of the condition
for the existence of an almost complex structure adapted to a surface as given in Lemma 1 of \cite{B1}, namely:

\begin{kor}\label{cor2}
Let $(X,J)$ be an almost complex 4-manifold 
and assume we are given an embedded surface $F \subset X$. 
Then the following
conditions are equivalent:
\begin{enumerate}
\item There is a (generic) $J'$ homotopic to $J$ such that the algebraic
number of negative complex points of $F$ with respect to $J'$ is zero.
\item There is a $J''$ homotopic to $J$ such that $F$ is pseudoholomorphic
with respect to $J''$.
\end{enumerate}
\end{kor}

\begin{proof}
First suppose that the algebraic number $n^-$ of negative complex points  with respect to
$J'$ is zero. Then, according to equation~\eqref{laimod}, the adjunction equality is fulfilled
and, by Lemma 1 in \cite{B1}, we can find an almost complex structure $J''$ homotopic
to $J'$ and hence to $J$ such that $F$ is pseudoholomorphic with respect to $J''$.
If conversely $F$ is pseudoholomorphic with respect to $J''$, the adjunction equality
holds for $J''$, and if we choose a generic $J'$ homotopic to $J$ (and $J''$), 
equation~\eqref{laimod}
implies that $n^-=0$.
\end{proof}

At this point, the author would like to point out that the proof of Lemma 2 in \cite{B1} contains
a minor gap, the arguments given there do not work in the special case $b^+=2=b^-$. However,
the assertion of the Lemma is true and a slight modification of the published version of the proof 
also works in this special case. A corrected version can be found in \cite{B2}.

Now let us relate Lai's results to the double points of immersed surfaces.
In the sequel, we will always assume that an immersion of a surface is
proper in the sense that the only singularities are ordinary double points.
Recall that there is a natural way to attach a sign to a double point $p$
of an oriented surface, depending on whether the orientations of the two
branches of the surface meeting in $p$ fit together to give the orientation
of $X$ at this point (then the sign should be $+1$) or not (sign $-1$). 
If the sign of a double point is $+1$,
we will call this point a {\bf positive double point}, otherwise it will
be called a {\bf negative double point}. We will need the following relation between
the double points of an immersed surface, its normal Euler number and its self--intersection number.

\begin{lemma}\label{lemma1}
Assume that $X$ is a closed and oriented
4--manifold and $\iota : F \hookrightarrow  X$ an immersion 
of a connected and oriented surface
having $d_+$ positive and $d_-$ negative double 
points. Let $N \rightarrow F$ denote 
the normal bundle of the immersion. Then
\[
e(N)=F \cdot F - 2 \, d_+ + 2 \, d_-.
\]
\end{lemma}

\begin{proof}
For the sake of simplicity, let us assume that there is only one double point $p$.
First consider the case that $p$ is positive. Let $\iota^{-1}(p)=\{x_1,x_2\}$
and choose an orientation preserving chart $h : U \rightarrow \R^4$ around $p$
such that
\begin{enumerate}
\item $h(U \cap \iota(F))=h(U) \cap (\R^2 \times 0 \cup 0 \times \R^2)$
\item For small disks $D_i$ in $F$ around $x_i$, $h(\iota(D_1)) \subset \R^2
\times 0$,$h(\iota(D_2)) \subset 0 \times \R^2$ and the restrictions $(h \circ \iota) | _{D_i}$
map the orientation of $F$ to the canonical orientations of the planes
$\R^2 \times 0$ and $0 \times \R^2$.
\end{enumerate} 
Then a trivialization for the normal bundle $N=\iota^*TM/TF$ restricted to 
$D_1$ is given by $dh \circ d\iota$, followed by the projection onto the 
second plane $0 \times \R^2$, and with respect to this chart, a section of
$N | _{D_1}$ is given by the affine plane $\R^2 \times 0 + \epsilon$ for
a small $\epsilon > 0$. If we choose a tubular neighborhood $\tau : N \hookrightarrow X$ 
which coincides with the map given by $dh$ and $h^{-1}$ around the $x_i$,
the image of this section will be contained in 
$h^{-1}(\R^2 \times 0 + \epsilon)$ and will intersect $F$ in one positive 
point. A similar section can be constructed over $D_2$, and combined with a
generic section of $N$ outside of the $D_i$, we obtain an immersion $\iota'$
of $F$ which will intersect $\iota(F)$ in $2+e(N)$ points, counted with signs.

If the double point $p$ is negative, the sections constructed over the $D_i$
will contribute with sign $-1$ to the intersection number of $\iota'(F)$
and $\iota(F)$, and we obtain $F \cdot F = e(N)-2$. This proves our
assertion in the case that there is only one double point, the proof in the
general case is similar.
\end{proof}

Now let us combine Lai's work -- namely equation~\eqref{laiorg} -- 
with Lemma~\ref{lemma1} to derive a relation between the homology class
of an immersed surfaces, its genus, the number of complex points
and the number of double points. We then obtain the following

\begin{proposition}\label{base} 
Let $(X,J)$ be an almost complex 4--manifold and $F \hookrightarrow X$ a
generic immersion of an oriented surface $F$ 
having $d_+$ positive and $d_-$ negative double points. Let
$n^-$ denote the algebraic number of negative complex points (as defined by
equation~\eqref{defnofalgsum}). Then
\[
g(F)+n^- - d_- + d_+ = 1 + \frac{1}{2} (F \cdot F + K \cdot F).
\]
\end{proposition}

\begin{proof}
We have a decomposition
$\iota^*TX = TF \oplus N$,
to which we can apply equation~\eqref{laiorg} to obtain
\[
(2-2g) - c_1(J) \cdot F + e(N) = 2 n^-.
\]
By Lemma~\ref{lemma1}, $e(N)=F \cdot F - 2d^+ + 2d^-$. Substituting this
into the last equation leads to the desired result. 
\end{proof}

\begin{rem}\label{changingor}
Assume that $F$ is an immersed surface with $d_{\pm}(F)$ double points
and $n^{\pm}(F)$ complex points. Let $\bar{F}$ denote the same surface with the reversed
orientation. Then clearly $d_{\pm}(\bar{F})=d_{\pm}(F)$, and hence Proposition~\ref{base}
yields the relations $n^+(\bar{F})=n^-(F)$ and $n^{-}(\bar{F})=n^+(F)$.
\end{rem}

\section{Immersed surfaces and genus estimates}\label{estimates}

Proposition~\ref{base} can be used to derive estimates for the number of negative complex points of immersed
surfaces (note that these estimates always include estimates on the algebraic number of
positive complex points, since Theorem~\ref{laithm} and equation~\eqref{laiorg} together imply
$n^- - n^+ = K \cdot F$).
As a first example, we will now
prove a bound for the number of double points of an immersion by
using branched covers as in~\cite{Ro}.

\begin{proposition}\label{prop1}
Let $X$ be a closed oriented and simply--connected
4--manifold and $F \subset X$ be an immersed surface of genus g
having $d_+$ positive
and $d_-$ negative double points. Assume that the homology class $[F]$ of
$F$ is divisible by a prime power $m$. 
Then
\[
d_++g \geq \frac{m+1}{6m}F \cdot F - b_2^+(X).
\]
\end{proposition}

\begin{proof}
As demonstrated in \cite{FS}, we can
find an immersed surface in $X \# \cpq$ having $d_+$ positive double
points, $d_- -1$ negative double points and representing the class 
$([F],0) \in H_2(X \# \cpq;\Z)= H_2(X;\Z)\oplus \Z$. In fact, pick two generic lines
in $\cpq$ and reverse the orientation of one of them to obtain two spheres $S_1, S_1$ which intersect
in one point with intersection number $+1$. Remove a small ball $B$ around this point and
a similar ball $B'$ around one negative double points of $F$. We than can glue $X \setminus B'$
and $\cpq \setminus B$ along their boundaries in such a way that the boundary links $F \cap B'$ and
$ (S_1 \cup S_2 ) \cap B$ get identified. This yields a new immersed surface in $X \# \cpq$ as desired.
By iterating the construction, one can construct an immersed surface in $X'=X \# d_- \cpq$ having
$d_+$ positive self--intersection points and representing the class $([F],0, \dots, 0)$.

Cutting out small disks and gluing in handles at the remaining 
double points leads to a surface $\Sigma \subset X'$ having genus $d_++g$
and self--intersection number $F \cdot F$. Since the homology class $[F]$ was divisible
by $m$, the same is true for $[\Sigma]$. Now let $Y \rightarrow X'$
denote the branched cover of order $m$ with branch locus $\Sigma$. As in 
\cite{Ro}, one can use this cover to obtain estimates for the genus
of $\Sigma$.  
We have
\begin{align*}
b_2(Y)&=m b_2(X') +2(m-1)(d_++g) \\ &= m b_2(X)+m d_-+2(m-1)(d_++g) 
\end{align*}
and 
\[
\tau(Y)=m\tau(X)-md_- - \frac{m^2-1}{3m} F \cdot F.
\]
Now the real cohomology $H^2(X;\R)$ appears in $H^2(Y;\R)$ as the subspace 
invariant under the action of $\Z_m$, and the splitting of $H^2(Y;\R)$ into
the eigenspaces of the action is orthogonal with respect to the intersection
form (see \cite{Ro}), hence we have the inequality
\[
|\tau(Y)-\tau(X)| \leq b_2(Y)-b_2(X),
\] 
in particular
\[
-\tau(Y)+\tau(X) \leq b_2(Y)-b_2(X).
\]
If we substitute the values for $b_2$ and $\tau$ from the above equation, we see
that the terms containing $d_-$ cancel out,
and this leads to the desired result.
\end{proof}

Note that reversing the orientation of $X$ gives an estimate for the number
of negative double points.
Furthermore,
the estimate of Rohklin (inequality 6.3. in \cite{Ro}) appears in the result as the
special case $d_+=d_-=0$.

\begin{ex} There is a simple example where the estimate in Proposition~\ref{prop1}
is sharp.
Consider the class $a=(3,3) \in H_2(\CP^2\#\CP^2;\Z)$. 
Suppose there is an
immersed sphere $S$ representing $a$ which has $d_+$ positive and
$d_-$ negative double points. 
Proposition~\ref{prop1} then shows that $d_+ \geq 2$.
On the other hand, the class $3[\CP^1] \in H_2(\CP^2;\Z)$ can be realized by a sphere with
one positive self--intersection point.
Taking two copies of this immersed sphere shows that the class $(3,3)$ can be represented
by a sphere having two positive double points. Hence the bound from Proposition
\ref{prop1} is sharp in this case. 
\end{ex}

\begin{ex}\label{kmexplained}
The result of Fintushel and Stern in \cite{FS} can be summarized by the 
statement that
for a sphere immersed in a rational surface, the minimal number of positive
double points is at least the minimal genus of its homology class. The estimate
in Proposition~\ref{prop1} goes in the same direction, since its right hand side
equals the right hand side in \cite{Ro}, 6.3 (up to the absolute value). There is a 
simple example showing that in general, one cannot estimate the number
of positive double points by the minimal genus (the same is true for the
number of negative double points).

First note that the  minimal genus of the class $a=(3,1) \in H_2(\CP^2 \# \CP^2;\Z)$ is one.
In fact, a torus representing $a$ is obtained by tubing together algebraic representatives
in both factors,
and by \cite{KeM}, $a$ cannot be represented by a sphere.
In the same paper,  Kervaire and Milnor show how to represent the class $(3,0)$ by
a sphere. From a configuration of three lines in $\CP^2$, we obtain a
sphere with one positive double point representing $3[\CP^1] \in H_2(\CP^2;\Z)$.
Take two lines in $\CP^2$ and reverse the orientation of one of them 
to produce two spheres that intersect each other in one negative intersection
point. Cutting out balls around the intersection points and identifying
the boundaries leads to a sphere $S \in \CP^2 \# \CP^2$ that represents
$(3,0)$.

From the construction one can deduce that $S$ will intersect a generic
line $\gamma$ in the second factor twice in one negative and one positive
intersection point. Tubing together at the positive point leads to an immersed
sphere representing the class $a$ with $d_+=0,d_-=1$.
This provides an example of an immersed sphere for which the number of
positive self--intersection points is smaller than the minimal genus in
this homology class. Gluing at the other
point gives a representative of the same class with $d_+=1,d_-=0$,
thus showing that also the number of negative intersection points can be smaller
than the minimal genus.
Reversing the orientation of the line finally leads to immersed spheres representing
the class $(3,-1)$.
\end{ex}

\smallskip

As was already indicated above, we can now combine the estimate given in Proposition
\ref{prop1} and the relation between complex points and double points of an immersion
to obtain the following result.

\begin{kor}\label{korwithro}
Suppose that $X$ is a simply--connected almost complex 4--mani\-fold with canonical 
class $K$ and $F \hookrightarrow X$ an immersion of a surface of genus $g$
having $d_+$ positive and $d_-$ negative
double points. Assume that the homology class of $F$ is divisible by a prime
power $m$. As usual, let $n^-$ denote the algebraic number of negative complex
points. Then
\[
d_- - n^- \geq \frac{1-2m}{6m} F \cdot F - \frac{1}{2} K \cdot F - b_2^+(X)-1.
\]
\end{kor}

\begin{proof}
This follows immediately from Proposition~\ref{prop1} and Proposition~\ref{base}.
\end{proof}

\begin{ex}
Consider the homology class $3[\CP^1] \in H_2(\C P^2;\Z)$. Suppose that $S$ is an immersed
sphere representing this class. Then $S$ cannot be an embedding, because the minimal
genus in this homology class is $1$ (see \cite{KeM}), hence there must be at least
one double point (as always, we assume that the immersion is generic with double points
as the only singularities). In fact, Proposition~\ref{prop1} implies that at least one
of the double points must have sign $+1$ (this follows also from \cite{FS}), i.e. $d_+ > 0$.
By Corollary~\ref{korwithro}, we have $n^- \leq d_-$. 

Since the estimate of Rokhlin is sharp for this homology class, we expect that the
inequality $n^- \leq d_-$ is also sharp, which turns out to be true. 
For an example where equality occurs consider the
algebraic curve of degree $3$ given by the equation $x^3 + y^3 = 3xyz$. This curve has one
node at the point $[0:0:1]$. From our point of view, it defines an immersion of a sphere
with one positive double point which  represents the homology class $3[\CP^1] \in H_2(\C P^2;\Z)$
and is pseudoholomorphic with respect to the canonical almost complex structure on $\C P^2$.
Hence we have $n^-=d_-=0$.
\end{ex}

\section{Proof of Theorem \ref{diskbundle}}\label{proof}

It is sufficient to prove the first assertion, the second then
follows by reversing the orientation of $F$, using Remark~\ref{changingor}. 
Let $T$ denote a tubular neighborhood of $F_0$ in which $F$ is contained.
Since $H_2(T;\Z)=H_2(F_0;\Z)=\Z$, there is a unique integer $k$ such that
$k[F_0]=[F] \in H_2(T;\Z)$. Clearly this relation also holds in $H_2(X;\Z)$
and therefore $k > 0$ since $F_0 \cdot F = k F_0 \cdot F_0 > 0$ and
$F_0 \cdot F_0 > 0$ by assumption.
First we will prove the estimate
\begin{equation}\label{localthomconj}
d_+(F) + g(F) \geq 1 + k (g(F_0)-1) + \begin{pmatrix} k \\ 2 \end{pmatrix}
F_0 \cdot F_0.
\end{equation}
Note that this is a version of the so-called ``local Thom conjecture'' for
immersed surfaces. In the special case of embedded surfaces, this conjecture
has been proved by Lawson (see \cite{La}). 

For the proof of \eqref{localthomconj}, we will not make use of the
almost complex structure on $X$ and of the fact that $F_0$ is a pseudoholomorphic
curve.
By general position, we can assume 
that $F$ and $F_0$ intersect transversely,
in particular no double point of $F$ is lying on $F_0$.
Pick a complex structure on $F_0$ and a holomorphic line
bundle $L$ over $F_0$ having degree $\degree(L)=\langle c_1(L) , [F_0] \rangle = F_0 \cdot F_0$. 
If we choose a metric on $L$, we can identify $T$ with the unit disk bundle of $L$.
Let $E=L \oplus \C$ and $Y={\mathbb P}E$ the total space
of the projective bundle associated to $E$. Then $Y$ is an algebraic
surface, and $b_2^+(Y)=1+2p_g(Y)=1$ (see for instance\cite{BPV} IV.2.6)

We have an embedding $L \hookrightarrow Y$, given by $l \mapsto [(l,1)]$,
and the image of $F_0$ under this embedding -- that again will be
denoted by $F_0$ -- is an algebraic curve in $Y$ having self--intersection
number $F_0 \cdot F_0=\degree(L)$. If $K$ denotes the canonical class
of $Y$, this implies $K \cdot F_0 = 2g(F_0)-2-F_0 \cdot F_0$.

Now let $Y'=Y \# d_-(F)\cpq$, where the blow--up is performed at the
positive self--intersection points of $F$. 
Using the construction of \cite{FS} as in the proof of Proposition~\ref{prop1},
we can construct an embedded surface in $Y'$
representing the homology class $[F]=k[F_0]$ having $d_+$ positive
double points (here we think of $H_2(Y;\Z)$ as a subgroup of $H_2(Y';\Z)$).
Replacing the remaining double points by handles
leads to an embedded surface $F'$ of
genus $g(F)+d_+$ with $[F']=[F]$.
Note that 
$Y'$ is again an algebraic surface with canonical class $K'=K-\sum_i E_i$,
where $E_i$ denotes the exceptional curve in $\cpq$.

Choose a K\"ahler metric $g'$ on $Y'$ and let $\omega'$ denote its
fundamental form. Since $F_0 \subset Y'$ is holomorphic, we have $[\omega'] \cdot [F_0] > 0$.
Now the K\"ahler metric
$g'$ defines a symplectic structure with symplectic form $\omega'$,
and $[F'][\omega']=k[F_0][\omega'] > 0$, $F' \cdot F' =k^2 F_0 \cdot F_0 > 0$, 
hence we can apply Theorem E in \cite{LL} to conclude that
\[
g(F') \geq 1 + \frac{1}{2}(K' \cdot F' + F' \cdot F').
\]
Substituting the values for $g(F')$ and $[F']$ leads to
\[
d_+ + g(F) \geq 1 + \frac{1}{2}( k(2g(F_0)-2-F_0 \cdot F_0) +
k^2 F_0 \cdot F_0),
\]
and the inequality~\eqref{localthomconj} follows.

Let us now proceed with the proof of Theorem~\ref{diskbundle}.
The inequality~\eqref{localthomconj} is an estimate for the genus of $F$ and the number
of double points in terms of
its homology class. In a second step, we can now use Proposition~\ref{base} to obtain
an inequality between the number of complex points and the number of double points.
In fact, Proposition~\ref{base} and the inequality~\eqref{localthomconj} together imply
\[
d_-(F)-n^-(F) \geq k \, g(F_0) - \frac{1}{2} k \, 
K \cdot F_0 - k - \frac{1}{2} k \, F_0 \cdot F_0.
\]
The right hand side of this estimate can still be simplified using the fact 
that $F_0$ was assumed to
be pseudoholomorphic and thus the adjunction equality holds for it. 
Therefore we finally obtain
\[
d_-(F)-n^-(F) \geq 0,
\]
and this is the desired estimate.

\end{document}